# On a theorem of Rickards

Håvard Damm-Johnsen[1]

A recent result of Rickards states that the generating series of intersection numbers of real quadratic geodesics on indefinite Shimura curves are elliptic modular forms. We reinterpret this as a Kudla–Millson theta series, and prove that Rickards' generating series is the diagonal restriction of a Hilbert modular form, analogous to results of Darmon–Pozzi–Vonk and Branchereau.

## Contents

## 1 Introduction

Let $B$ be a non-split rational indefinite quaternion algebra of discriminant $D_B$, and let $\mathcal{O}$ be an Eichler order of level $N$. Fix an isomorphism $B \otimes \mathbb{R} \cong \mathrm{GL}_2(\mathbb{R})$. This identifies $\mathcal{O}^1 := \{x \in \mathcal{O} : \mathrm{Nm}(x) = 1\}$ with a discrete subgroup of $\mathrm{SL}_2(\mathbb{R})$, giving rise to a compact Shimura curve $X_{\mathcal{O}} := \mathcal{O}^1 \backslash \mathfrak{h}$, where $\mathfrak{h}$ denotes the complex upper-half plane.

Fix a pair of distinct algebra embeddings of real quadratic fields $\alpha_j : F_j \to B$, $j \in \{1, 2\}$. These naturally give rise to geodesics $\tilde{C}_{\alpha_j}$ in $\mathfrak{h}$, as explained in Section 2.1. If $\alpha_j(\mathcal{O}_{F_j}) \subset \mathcal{O}$, then a segment of the geodesic $\tilde{C}_{\alpha_j}$ naturally maps onto a closed cycle in $X_{\mathcal{O}}$, which we denote by $C_{\alpha_j}$.

The arithmetic interest of real quadratic geodesics in the upper half plane was recognized by Hecke in on his eponymous integral formula [Hec17], and they feature prominently in the Shintani's work on the theta correspondence [Shi75]. More recently, in the work of Darmon and Vonk on rigid meromorphic cocycles, such geodesics play a key role, in some respects analogous to that of CM points in CM theory. In their work, height pairings of CM points have a $p$-adic counterpart in certain $p$-adic height pairings of such geodesics. These may be thought of as a higher analogue of the oriented intersection product $\langle -, - \rangle_{X_\mathcal{O}}$ on $H_1(X_\mathcal{O}, \mathbb{Z})$, defined in terms of the cup product of the respective fundamental classes, see [GH81, §31]. Let $T_n$ denote the $n$-th Hecke correspondence acting on this homology group.[2]

---

[1]University of Gothenburg, Sweden

In [Ric22], Rickards proved the following theorem:

**Theorem 1.1** ([Ric22]): *There exists a modular form $R(\tau) = \sum_{n=1}^\infty a_n(R) q^n$ in $S_2(\Gamma_0(D_B N))$ with*

$$a_n(R) = \langle C_{\alpha_1}, T_n C_{\alpha_2} \rangle_{X_\mathcal{O}}, \tag{1.1}$$

*for all $(n, N) = 1$.*

The idea behind Rickards' proof is straightforward: the map $T_n \mapsto \langle C_{\alpha_1}, T_n C_{\alpha_2} \rangle$ defines a linear functional on the Hecke algebra acting faithfully on $H_1(X_\mathcal{O}, \mathbb{Z})$, and by Jacquet–Langlands, this transfers to a modular form with the corresponding Fourier coefficients. An interesting feature of Rickards' work is that the Hecke operators may in fact be defined directly on the set of optimal embeddings, enabling explicit computation.

An analogous generating series on the modular curve $Y_0(p)$ is considered in [DPV21], where the geodesic $\tilde{C}_{\alpha_1}$ is replaced with a segment of the line between $0$ and $\infty$ in $\mathfrak{h}$. They show (Theorem A) that the generating series is in fact the restriction to the diagonal of a certain Hilbert Eisenstein series associated with $F_2$. Their argument is essentially an explicit computation of Fourier coefficients. This was reinterpreted as a (regularized) Kudla–Millson theta lift in [Bra23a]. Our main result is that a similar statement holds for Rickards' generating series.

**Theorem 1.2**: *Assume that the discriminants $D_j$ of $F_j$ satisfy $(D_1, D_2) = 1$, and let $F := \mathbb{Q}\left(\sqrt{D_1 D_2}\right)$. Then there exists an explicit Hilbert modular form $\vartheta$ over $F$ of parallel weight 1 such that $R(\tau) = \vartheta(\tau, \tau)$.*

We construct $\vartheta$ explicitly, and consequently get an expression for the intersection numbers in terms of its Fourier coefficients. A key result for this part of the argument is the existence of an $F$-valued quadratic form on $B$, defined in terms of the embeddings $\alpha_1$ and $\alpha_2$, and which satisfies $\mathrm{Tr}_{F/\mathbb{Q}}\, q_F(b) = \mathrm{Nm}(b)$ for all $b \in B$. A little more work gives an isometry of $F$-quadratic spaces $\iota_L : B \to L$, where $L$ has the quadratic form $q(x) := q_F(1) \cdot \mathrm{Nm}_{L/F}(b)$.

---

[2]More precisely, this is the Hecke action compatible with the usual action on cohomology via Poincaré duality.

In the special case where $q_F(1)$ is totally positive, the formula for the Fourier coefficients of $R(\tau)$ reduces to the following:

**Corollary 1.3**: *Let $\varepsilon_1$ and $\varepsilon_2$ denote the fundamental totally positive units in $F_1$ and $F_2$ respectively, and suppose $q_F(1) \gg 0$. Then*

$$\langle C_{\alpha_1}, T_n C_{\alpha_2} \rangle_{X_{\mathcal{O}}} = \frac{1}{2} \sum_{\substack{b \in \varepsilon_1^{\mathbb{Z}} \backslash \mathcal{O} / \varepsilon_2^{\mathbb{Z}} \\ q_F(b) \gg 0 \\ \mathrm{Tr}_{F/\mathbb{Q}}\, q_F(b) = n}} \mathrm{sgn}\, \mathrm{Nm}_{F/\mathbb{Q}}\, \mathrm{Tr}_{L/F}(\iota_L(b)). \tag{1.2}$$

We now briefly sketch the main ideas of the proof. By [Lemma 2.8](), the intersection number $\langle C_{\alpha_1}, T_n C_{\alpha_2} \rangle$ equals $\langle C_{\alpha_1} \times C_{\alpha_2}, C_n \rangle$ in $X_{\mathcal{O}} \times X_{\mathcal{O}}$, where $C_n$ is graph of the $n$-th Hecke correspondence. Viewing $B$ as a rational quadratic space with the norm form, the locally symmetric space of $\mathrm{SO}_B$ may be identified with $X_{\mathcal{O}} \times X_{\mathcal{O}}$, and $C_n$ agrees with the $n$-th special cycle, in the sense of Kudla and Millson. The Kudla–Millson geometric theta kernel $\Theta_{\mathrm{KM}}$ is a modular form for $\mathrm{SL}_2$ valued in the cohomology of $X_{\mathcal{O}} \times X_{\mathcal{O}}$, and is characterized by its $q$-expansion

$$\Theta_{\mathrm{KM}}(\tau) = \sum_{n=1}^{\infty} \mathrm{PD}(C_n) q^n, \tag{1.3}$$

where $\mathrm{PD}(C_n)$ denotes the Poincaré dual of $C_n$. The cap product with $C_{\alpha_1} \times C_{\alpha_2}$ then gives Rickards' generating series $R$.

To show that $R$ is a diagonal restriction, an argument in [Section 3.1]() inspired by Howard and Yang [HY12] gives the isometry of $F$-quadratic spaces $\iota_L$ between $B$ and $L := F_1 \otimes_{\mathbb{Q}} F_2$ with a certain quadratic form. We reinterpret $C_{\alpha_1} \times C_{\alpha_2}$ as a torus in $\mathrm{SO}_B$ whose image in $\mathrm{SO}_L$ is identified with $\mathrm{Res}_{\mathbb{Q}}^F \mathrm{SO}_L$. Using this, a geometric theta lift from $\mathrm{Res}_{\mathbb{Q}}^F \mathrm{SO}_L$ to $\mathrm{SL}_2(F)$ gives a Hilbert modular form $\vartheta$, whose diagonal restriction equals $R$ by the following seesaw diagram:

$$\begin{array}{ccc} \mathrm{SO}_L & & \mathrm{Res}_{\mathbb{Q}}^F \mathrm{SL}_2 \\ \uparrow & \diagdown\mkern-14mu\diagup & \uparrow \\ \mathrm{Res}_{\mathbb{Q}}^F \mathrm{SO}_L & & \mathrm{SL}_2 \end{array}$$

This is a special case of the seesaw diagram in [Kud84, Equation (2.19)], to which the reader is referred for further details.

In the recent paper [Bra25], Branchereau extends his aforementioned result to a setting of general étale algebras. The proof of Theorem 1.2 essentially combines ideas of that paper and [Bra23b], with a couple of new ideas:

- We use ideas of Howard and Yang [HY12] from the CM case to get an explicit description of the étale algebra corresponding to our geodesics in terms of optimal embeddings, analogous to Rickards' setup. This in turn gives the Fourier coefficient expressions in Corollary 1.3.
- Our conceptual approach to Hecke correspondences in Section 2.2, particularly Lemma 2.8, replaces the direct computation in [Bra23b, §5.8].

Theorem 1.2 is in some sense a generalization of Theorem A in [DPV21], which is the case of the split quaternion algebra $B = \mathrm{Mat}_2(\mathbb{Q})$ and the "degenerate real quadratic field" $F_1 = \mathbb{Q} \times \mathbb{Q}$. This case was worked out and generalized in [Bra23a], albeit with a slightly different approach to the tori; cf. Remark 3.1(ii).

We hope that the theta lift perspective might shed light on possible generalizations of the other theorems in *loc. cit.* and their refinements in [DPV23]. In fact, this perspective already inspired the proof of an extension of the results of [DPV21] to ring class characters in [Dam24], though the proof there is somewhat ad-hoc.

**Acknowledgements.** It is a pleasure to thank Mike Daas, Luis Garcia, Alex Horawa, James Newton and Jan Vonk for insightful discussions, as well as Romain Branchereau for helpful email correspondence clarifying technical points regarding his work. At different stages of writing this paper the author was supported by an Aker scholarship and the Max Planck Institute for Mathematics in Bonn.

## 2 Generating series of geodesics

In this section, we prove that the Kudla–Millson theta lift recovers Rickards' Shimura curve generating series. We first introduce some conventions:

- For any number field $F$, $\mathbb{A}_F$ denotes the adele ring of $F$, and $\mathbb{A}_F^\infty$ the finite adeles. If $F = \mathbb{Q}$, we omit the subscript. Generally, $\infty$ as a subscript means "at $\infty$", while superscript $\infty$ means "away from $\infty$".
- If $G$ is an algebraic group over $F$, we set $[G] := G(F) \setminus G(\mathbb{A})$.
- Given a vector space $V$ over $F$ and an $F$-algebra $R$, we write $V_R$ for $V \otimes_F R$.
- If $k$ is a completion of $F$ or $\mathbb{A}_F$, we denote by $\mathcal{S}(V_k)$ the space of Schwartz–Bruhat functions on $V_k$.
- We fix the standard additive character $\psi : \mathbb{A} \to \mathbb{C}$ whose component at $\infty$ is $\psi_\infty(x) = e^{2\pi i x}$, and $\psi_p(x) = e^{-2\pi i \{x\}}$, where $\{\cdot\}$ is the fractional part. By construction, this factors through $\mathbb{Q} \setminus \mathbb{A}$. For any number field $F$, its precomposition with the trace map is an additive character $\psi_F : F \setminus \mathbb{A}_F \to \mathbb{C}$.

- If $(V, q)$ is a quadratic space, then we denote by

$$\langle v, w \rangle_V := q(x+y) - q(x) - q(y) \tag{2.4}$$

the associated bilinear form, so that $\langle v, v \rangle_V = 2q(v)$.

## 2.1 Symmetric spaces

We first recall some background on symmetric spaces for orthogonal groups, in the special case of the quadratic space $B$ with the norm form $\mathrm{Nm} : B \to \mathbb{Q}$. For a general overview, see for example [Kud03, §1]. Note that $(B, \mathrm{Nm})$ is an anisotropic four-dimensional rational quadratic space with real signature $(2,2)$. Set $G := \mathrm{SO}_B$, and let $H := \mathrm{GSpin}_B$, which by definition fits into a short exact sequence of algebraic groups,

$$1 \to \mathbb{G}_m \to H = \mathrm{GSpin}_B \to \mathrm{SO}_B \to 1. \tag{2.1}$$

Explicitly, we may identify $H(\mathbb{Q})$ with

$$B^\times \times_{\mathrm{Nm}} B^\times := \{(b_1, b_2) \in B^\times \times B^\times : \mathrm{Nm}(b_1) = \mathrm{Nm}(b_2)\} \tag{2.2}$$

with the map to $\mathrm{SO}_B(\mathbb{Q})$ coming from the action on $B$ given by $(h_1, h_2) \cdot b = h_1 b h_2^{-1}$. Let $\mathbb{D}$ denote the real symmetric space of $\mathrm{SO}_B(\mathbb{R})$, consisting of oriented negative 2-dimensional subspaces of $B_\mathbb{R}$. Under the identification $B_\mathbb{R} \cong \mathrm{Mat}_2(\mathbb{R})$, the norm on $B$ is identified with the determinant, and we fix the base point

$$X_0 := \left\{ \begin{pmatrix} x & y \\ y & -x \end{pmatrix} \subset \mathrm{Mat}_2(\mathbb{R}) : x, y \in \mathbb{R} \right\} \in \mathbb{D}. \tag{2.3}$$

Then $\mathrm{SO}_B(\mathbb{R})/\mathrm{Stab}(X_0)$ is isomorphic to the connected component of $X_0$ in $\mathbb{D}$, denoted $\mathbb{D}^+$. The stabilizer of $X_0$ in $H(\mathbb{R})^+ \cong \mathrm{GL}_2(\mathbb{R})^+ \times_{\det} \mathrm{GL}_2(\mathbb{R})^+$ is $\mathbb{R}_{>0} K_{H,\infty}$ where $K_{H,\infty} := \mathrm{SO}(2) \times \mathrm{SO}(2)$, and so $\mathbb{D}^+$ is naturally the symmetric space of $H$. On the other hand, $\mathrm{GSpin}_B(\mathbb{R})^+$ acts naturally on $\mathfrak{h} \times \mathfrak{h}$ by Möbius transformations in each variable, and the point $(i, i)$ has stabilizer $K_{H,\infty}$. This gives an identification $\mathbb{D}^+ \cong \mathfrak{h} \times \mathfrak{h}$. Explicitly, a plane $z = h_1 X_0 h_2^{-1} \in \mathbb{D}^+$ is mapped to $(h_1 i, h_2 i)$.

Fix an Eichler order $\mathcal{O}$ in $B$ and a pair of algebra embeddings $\alpha_1 : F_1 \to B$ and $\alpha_2 : F_2 \to B$, where $F_1$ and $F_2$ are real quadratic fields with coprime discriminants. Then each $\alpha_j$ is an optimal embedding of the order $\mathcal{O}_j := \alpha_j^{-1}(\mathcal{O})$ in $F_j$. This gives rise to a torus $T$ defined by

$$T := \mathrm{Res}_\mathbb{Q}^{F_1} \mathbb{G}_m \times_{\mathrm{Nm}} \mathrm{Res}_\mathbb{Q}^{F_2} \mathbb{G}_m. \tag{2.4}$$

More precisely, this is the algebraic group which has $R$-points

$$T(R) = \{(t_1, t_2) \in (F_1 \otimes R)^\times \times (F_2 \otimes R)^\times : \mathrm{Nm}_1(t_1) = \mathrm{Nm}_2(t_2)\}, \quad (2.5)$$

for any $\mathbb{Q}$-algebra $R$. Here $\mathrm{Nm}_j : F_j \otimes R \to R$ is the corresponding norm map, which agrees with the usual norm on $F_j$ when $R = \mathbb{Q}$. The torus $T$ is the real quadratic analogue of the torus $T$ considered in [HY12], and there is a natural embedding of $T$ into $H$ given by $(t_1, t_2) \mapsto (\alpha_1(t_1), \alpha_2(t_2))$.

To find the symmetric space of $T$, denoted $\mathbb{D}_T$, pick isomorphisms $\gamma_j : F_j \otimes \mathbb{R} \cong \mathbb{R} \times \mathbb{R}$. We identify $\gamma_j$ with the corresponding matrix in $\mathrm{GL}_2(\mathbb{R})$, and let $\gamma_\infty := (\gamma_1, \gamma_2)$. Then $T(\mathbb{R}) = \gamma_\infty T_\Delta(\mathbb{R}) \gamma_\infty^{-1}$, where

$$T_\Delta(\mathbb{R}) := \{((x_1, y_1), (x_2, y_2)) \in (\mathbb{R}^\times \times \mathbb{R}^\times)^2 : x_1 y_1 = x_2 y_2\}. \quad (2.6)$$

As a subgroup of $\mathrm{GL}_2(\mathbb{R}) \times_{\det} \mathrm{GL}_2(\mathbb{R})$, $T_\Delta(\mathbb{R})$ is the preimage of the diagonal torus in $\mathrm{GL}_2(\mathbb{R}) \times \mathrm{GL}_2(\mathbb{R})$, and its maximal compact subgroup is $K_{T_\Delta, \infty} := \{\pm 1, \pm 1\}$, where the signs are independent. It follows that

$$\mathbb{D}_T \cong \mathbb{D}_{T_\Delta} := T_\Delta(\mathbb{R}) / \left(\mathbb{R}^\times K_{T_\Delta, \infty}\right) \cong (\mathbb{R}_{>0})^2. \quad (2.7)$$

We view $\mathbb{D}_T$ as a subspace of $\mathbb{D}$ by picking the base point $z_0 = \gamma_\infty \cdot X_0$, noting that $K_{T, \infty} := \gamma_\infty K_{T_\Delta, \infty} \gamma_\infty^{-1}$ stabilizes $z_0$.

**Proposition 2.1**: *For $j = 1, 2$, let $\tilde{C}_{\alpha_j}$ be the geodesic in $\mathfrak{h}$ connecting $\gamma_j(0)$ and $\gamma_j(\infty)$. Then the image of $\mathbb{D}_T$ in $\mathfrak{h} \times \mathfrak{h}$ is $\tilde{C}_{\alpha_1} \times \tilde{C}_{\alpha_2}$.*

*Proof*: Fix $t_\Delta = ((t_1, t_1^{-1}), (t_2, t_2^{-1})) \in T_\Delta(\mathbb{R})$, and let $t := \gamma_\infty t_\Delta \gamma_\infty^{-1}$ be the corresponding element of $T(\mathbb{R})$. Then $t \cdot z_0 = \gamma_\infty t_\Delta \cdot X_0$, which maps to $(\gamma_1 t_1^2 i, \gamma_2 t_2^2 i)$ in $\mathfrak{h} \times \mathfrak{h}$. Let $I_{\{0, \infty\}}$ denote the geodesic in $\mathfrak{h}$ from $0$ to $\infty$. Varying $t_1$ and $t_2$ in $\mathbb{R}^\times$, the translate of $I_{\{0, \infty\}} \times I_{\{0, \infty\}}$ under $\gamma_\infty$ is precisely $\tilde{C}_{\alpha_1} \times \tilde{C}_{\alpha_2}$. □

We now introduce level structures. Let $\hat{\mathcal{O}} := \mathcal{O} \otimes_\mathbb{Z} \hat{\mathbb{Z}}$, and set

$$K_H^\infty := \left(\hat{\mathcal{O}} \times \hat{\mathcal{O}}\right) \cap \left((B \otimes \mathbb{A}^\infty)^\times \times_{\mathrm{Nm}} (B \otimes \mathbb{A}^\infty)^\times\right). \quad (2.8)$$

Then $K_H^\infty$ is a compact open subgroup of $H(\mathbb{A}^\infty)$, and by Eichler's theorem, $H$ has class number 1. Evidently, the image of $K_H^\infty$ in $\mathrm{SO}_B(\mathbb{A}^\infty)$ stabilizes the lattice $\hat{\mathcal{O}} \subset B \otimes \mathbb{A}^\infty$. It follows that the locally symmetric space for $H$ and $\mathrm{SO}_B$ is given by

$$X_H := \Gamma_H \backslash \mathbb{D} \quad \text{for} \quad \Gamma_H := (H(\mathbb{Q}) \cap K_H^\infty) \cap H(\mathbb{R})^+. \quad (2.9)$$

Via the isomorphism $\mathbb{D}^+ \cong \mathfrak{h} \times \mathfrak{h}$, it is easy to check that $X_H \cong X_\mathcal{O} \times X_\mathcal{O}$, where $X_\mathcal{O} = \mathcal{O}^1 \backslash \mathfrak{h}$. Similarly, we let $K_T^\infty$ be the preimage of $\hat{\mathcal{O}}_1 \times \hat{\mathcal{O}}_2$ in $T(\mathbb{A}^\infty)$, and define

$$\mathrm{Cl}\, T := T(\mathbb{Q}) \backslash T(\mathbb{A}^\infty)/K_T^\infty, \tag{2.10}$$

as well as

$$\mathrm{Cl}^+ T := T(\mathbb{Q}) \backslash T(\mathbb{A})/K_T^\infty T(\mathbb{R})^+, \tag{2.11}$$

where $T(\mathbb{R})^+ := T(\mathbb{R}) \cap (F_1 \otimes \mathbb{R})^+ \times (F_2 \otimes \mathbb{R})^+$.

**Lemma 2.2**: *There is an isomorphism*

$$T(\mathbb{Q}) \backslash T(\mathbb{A}^\infty)/K_T^\infty \to \mathrm{Cl}\, \mathcal{O}_1 \times \mathrm{Cl}\, \mathcal{O}_2. \tag{2.12}$$

*Proof*: The corresponding argument for imaginary quadratic fields given in [HY12, Prop. 2.14] applies with minor modifications. □

Let $X_T := T(\mathbb{Q}) \backslash T(\mathbb{A}^\infty) \times \mathbb{D}_T/K_T^\infty$ be the corresponding locally symmetric space. The group $\Gamma_T := (K_T^\infty \cap T(\mathbb{Q})) \cap T(\mathbb{R})^+$ acts discretely on $\mathbb{D}_T$. Explicitly, it is generated by $\varepsilon_1$ and $\varepsilon_2$, the totally positive fundamental units in $\mathcal{O}_1^\times$ and $\mathcal{O}_2^\times$. There is a diffeomorphism $X_T \cong \bigsqcup \Gamma_T \backslash \mathbb{D}_T$, obtained by sending $z$ in the component of $t^\infty$ to $(t^\infty, z)$ in $T(\mathbb{A}^\infty) \in \mathbb{D}_T$. The right-hand side is indexed by $\mathrm{Cl}^+ T$, or equivalently $\mathrm{Pic}^+ \mathcal{O}_1 \times \mathrm{Pic}^+ \mathcal{O}_2$, by Lemma 2.2. The following lemma shows that our original choice of optimal embeddings $\alpha_j$ is not essential.

**Corollary 2.3**: *The group $\mathrm{Cl}^+ T$ acts transitively on the set of pairs of oriented optimal embeddings $F_1 \times F_2 \to B \times B$.*

*Proof*: This follows from the corresponding statement about class groups, see for example [Ric20, Theorem 4.5.4], wherein the precise definition of "oriented" is also found. □

Since the inclusion $X_T \to X_H$ maps the connected component corresponding to $1 \in \mathrm{Cl}_T^+$ onto $\varepsilon_1^\mathbb{Z} \times \varepsilon_2^\mathbb{Z} \backslash \tilde{C}_{\alpha_1} \times \tilde{C}_{\alpha_2}$, we have shown:

**Corollary 2.4**: *The inclusion $X_T \to X_H \cong X_\mathcal{O} \times X_\mathcal{O}$ identifies the component of $1$ with the cycle $C_{\alpha_1} \times C_{\alpha_2} \in X_\mathcal{O} \times X_\mathcal{O}$.*

### 2.2 Special cycles and Hecke correspondences

The locally symmetric space $X_H$ has a natural supply of 2-cycles known as *special cycles*. It is well-known that these are embedded Shimura curves, see for example [DY13]; in this section we show that they may be identified with graphs of Hecke correspondences.

For $x \in B$ of positive norm, let

$$\mathbb{D}_x := \{z \in \mathbb{D}^+ : z \subset x^\perp\} \subset \mathbb{D}^+, \tag{2.1}$$

where $x^\perp$ is the orthogonal complement in $B(\mathbb{R})$ of the line spanned by $x$.

**Lemma 2.5**: *The image of $\mathbb{D}_x$ in $\mathfrak{h} \times \mathfrak{h}$ is the set $\{(x\tau, \tau) : \tau \in \mathfrak{h}\}$.*

*Proof*: A choice of orthogonal decomposition $B = \mathbb{Q}x \oplus x^\perp$ identifies the stabilizer $H_x := \mathrm{Stab}_{H(\mathbb{R})}(x)$ with $\mathrm{GSpin}_{x^\perp}$. Explicitly, since $g_1 x g_2^{-1} = x$ if and only if $g_1 = x g_2 x^{-1}$, this is given by

$$\mathrm{Stab}_{H(\mathbb{R})}(x) = \{(xgx^{-1}, g) : g \in B_\mathbb{R}\} = H_x(\mathbb{R}) \cong \mathrm{GL}_2(\mathbb{R}). \tag{2.2}$$

Under this isomorphism, we identify $\mathbb{D}_x$ with

$$\mathbb{D}_{H_x} = \{z \subset x^\perp : \dim z = 2 \text{ and } \mathrm{Nm}|_z < 0\} \tag{2.3}$$

by noting that the base point $xX_0$ is in $\mathbb{D}_x$. Indeed, multiplying by $x^{-1}$ reduces the problem to showing that $1 \in X_0^\perp$, which follows from a quick determinant computation. Then the map $\mathbb{D}_{H_x} \to \mathbb{D}_x$ is given by $g \cdot xX_0 \mapsto (xgx^{-1}, g) \cdot xX_0$. Furthermore, we identify $\mathbb{D}_{H_x}$ with $\mathfrak{h}$ as above, sending $\tau = g_\tau \cdot i \in \mathfrak{h}$ to $g_\tau \cdot xX_0$. The composite map

$$\mathfrak{h} \to \mathbb{D}_{H_x} \to \mathbb{D}_x \hookrightarrow \mathbb{D} \to \mathfrak{h} \times \mathfrak{h} \tag{2.4}$$

is given by

$$\tau \mapsto g_\tau \cdot xX_0 \mapsto (xg_\tau x^{-1}, g_\tau) \cdot xX_0 = (xg_\tau, g_\tau) \cdot X_0 \mapsto (x\tau, \tau), \tag{2.5}$$

which proves the claim. $\square$

Letting $\Gamma_x$ denote the stabilizer of $x$ in $\Gamma_H$, we get a submanifold $C_x := \Gamma_x \backslash \mathbb{D}_x$ of $X_H$. Define $\mathcal{O}_x := \mathcal{O} \cap x^{-1}\mathcal{O}x$. Then the map $\mathbb{D}_x \to \mathcal{O}_x^1 \backslash \mathfrak{h} =: X_{\mathcal{O}_x}$ given by $(x\tau, \tau) \mapsto \tau$ induces an isomorphism $C_x \cong X_{\mathcal{O}_x}$.

**Definition 2.6**: Fix $n \in \mathbb{Q}_{>0}$ and $\varphi^\infty \in \mathcal{S}(B_{\mathbb{A}^\infty})$. We define a **special cycle** as

$$C_n(\varphi^\infty) := \sum_{\substack{x \in \Gamma_H \backslash B \\ \mathrm{Nm}(x) = n}} \varphi^\infty(x) C_x. \tag{2.6}$$

Note that the a priori infinite sum is finite, so $C_n(\varphi^\infty)$ defines a genuine 2-cycle on $X_H$, which is a formal sum of embedded Shimura curves. Unfolding definitions, we find:

**Proposition 2.7**: *Let $\varphi^\infty = \mathbb{1}_{\hat{\mathcal{O}}}$, and let $C_n := C_n(\varphi^\infty)$. Then the image of $C_n$ in $X_\mathcal{O} \times X_\mathcal{O}$ is given by*

$$C_n = \sum_{x \in \mathcal{O}^n \,/\!/\, \mathcal{O}^1} C_x, \tag{2.7}$$

where $\mathcal{O}^n \,/\!/\, \mathcal{O}^1 := \{x \in \mathcal{O}^1 \setminus \mathcal{O}/\mathcal{O}^1 : \mathrm{Nm}(x) = n\}$.

In particular, $C_n$ is empty when $n = 0$ since $B$ is a division ring, and when $n$ is not an integer. Recall from [Zha01, §1.4] that the $n$-th Hecke correspondence $T_n$ is the sum over $x \in \mathcal{O}^n \,/\!/\, \mathcal{O}^1$ of geometric correspondences

$$\begin{array}{ccc} & X_{\mathcal{O}_x} & \\ {}^\alpha \swarrow & & \searrow {}^\beta \\ X_\mathcal{O} & & X_\mathcal{O} \end{array}$$

where $\alpha$ is the map induced by the inclusion $\mathcal{O}_x \to \mathcal{O}$, and $\beta$ comes from multiplication by $x$. The graph of $T_n$ is defined to be

$$\Gamma_{T_n} := \{(x_1, x_2) \in X_\mathcal{O} \times X_\mathcal{O} : x_2 \in \beta(\alpha^{-1}(x_1))\}, \tag{2.8}$$

and the two projection maps $p_1, p_2 : X_\mathcal{O} \times X_\mathcal{O} \to X_\mathcal{O}$ are covering maps with finite fibres. From the explicit description of $C_n$ under the identification $X_H \cong X_\mathcal{O} \times X_\mathcal{O}$, we see that $\Gamma_{T_n} = C_n$ pointwise.

**Lemma 2.8**: *Fix geodesics $C$ and $C'$ in $X_\mathcal{O}$. Then*

$$\langle C \times C', C_n \rangle = \langle C, T_n C' \rangle. \tag{2.9}$$

*Proof*: Let $H^\bullet_{\mathrm{dR}}(X_\mathcal{O})$ denote the de Rham cohomology ring of $X_\mathcal{O}$. We will show that

$$\int_{X_\mathcal{O} \times X_\mathcal{O}} (p_1^* \alpha \wedge p_2^* \beta) \wedge [C_n] = \int_{X_\mathcal{O}} \alpha \wedge T_n \beta, \tag{2.10}$$

for all $\alpha, \beta \in H^\bullet_{\mathrm{dR}}(X_\mathcal{O})$, where $[C_n] \in H^n_{\mathrm{dR}}(X_\mathcal{O} \times X_\mathcal{O})$ denotes the Poincaré dual of the fundamental class of $C_n$. Since the Hecke action on $H_\bullet(X_\mathcal{O})$ is compatible with Poincaré duality, the result then follows by taking $\alpha = [C]$ and $\beta = [C']$ in $H^1_{\mathrm{dR}(X_\mathcal{O})}$.

Let $\{\omega_i\}$ be a basis for $H^\bullet_{\mathrm{dR}}(X_\mathcal{O})$, and let $\{\tau_j\}$ denote the dual basis, satisfying

$$\int_{X_\mathcal{O}} \omega_i \wedge \tau_j = \begin{cases} 1 & \text{if } i = j, \\ 0 & \text{otherwise.} \end{cases} \tag{2.11}$$

By [Tu23, Proposition 3],

$$T_n \tau_j = \sum_k a_{kj} \tau_k \quad \text{where} \quad a_{kj} = \int_{C_n} p_1^* \omega_k \wedge p_2^* \tau_j. \tag{2.12}$$

It follows that

$$\int_{X_\mathcal{O}} \omega_i \wedge T_n \tau_j = a_{ij} = \int_{X_\mathcal{O} \times X_\mathcal{O}} (p_1^* \omega_i \wedge p_2^* \tau_j) \wedge [C_n], \tag{2.13}$$

and Equation (2.10) now follows by writing $\alpha = \sum_i \alpha_i \omega_i$ and $\beta = \sum_j \beta_j \tau_j$. □

### 2.3 The Kudla–Millson form

Let $k$ be a local field, and let $\psi : k \to \mathbb{C}^\times$ be an additive character. If $(V, q)$ is a quadratic space over $k$, then the *Weil representation* (or *oscillator representation*) for the pair $(\mathrm{SL}_2(k), O_V(k))$, has underlying space $\mathcal{S}(V)$ and satisfies the formulas

$$\begin{aligned}
\omega(h)\varphi(v) &= \varphi(h^{-1}v) \quad \text{for } h \in O_V, \\
\omega\begin{pmatrix} 1 & x \\ 0 & 1 \end{pmatrix}\varphi(v) &= \psi(xq(v))\varphi(v) \quad \text{for } x \in k, \\
\omega\begin{pmatrix} a & 0 \\ 0 & a^{-1} \end{pmatrix}\varphi(v) &= |a|^{\dim V/2}\varphi(av) \quad \text{for } a \in k^\times.
\end{aligned} \tag{2.1}$$

We denote by the same symbol the oscillator representation over the adeles of a number field.

The *Kudla–Millson form* $\varphi_{\mathrm{KM}} \in \mathcal{S}(B_\mathbb{R}) \otimes \Omega^2(X_H)$ is a rapidly decaying closed differential form, first constructed in [KM90]. Branchereau [Bra23b] recently gave a new construction of $\varphi_{\mathrm{KM}}$ using the Mathai–Quillen formalism.

Fix a Schwartz–Bruhat function $\varphi^\infty \in \mathcal{S}(B_{\mathbb{A}^\infty})$ which is invariant under the induced action of $K_H^\infty$. Then we may form the *Kudla–Millson theta series*

$$\Theta_{\mathrm{KM}}(\tau, \varphi^\infty) = \sum_{v \in B} \omega(g_\tau)\varphi^\infty \otimes \varphi_{\mathrm{KM}}(v) \in \Omega^2(X_H), \tag{2.2}$$

where for $\tau = x + iy$ we let

$$g_\tau = \begin{pmatrix} 1 & x \\ 0 & 1 \end{pmatrix}\begin{pmatrix} y^{\frac{1}{2}} & 0 \\ 0 & y^{-\frac{1}{2}} \end{pmatrix} \in \mathrm{SL}_2(\mathbb{R}), \tag{2.3}$$

so that $g_\tau \cdot i = \tau \in \mathfrak{h}$. This has the remarkable property that for any compact 2-cycle $C$ in $X_H$,

$$\int_C \Theta_{\mathrm{KM}}(\tau, \varphi^\infty) = \sum_{n \in \mathbb{Q}_{\geq 0}} \int_C \Theta_n(\varphi^\infty) q^n, \tag{2.4}$$

where $\Theta_n(\varphi^\infty)$ is the Poincaré dual of the cycle $C_n(\varphi^\infty)$.

**Proposition 2.9**: *Let $\varphi^\infty := \mathbb{1}_{\hat{\mathcal{O}}}$. Then*

$$\int_{C_{\alpha_1} \times C_{\alpha_2}} \Theta_{\mathrm{KM}}(\tau, \varphi^\infty) = \sum_{n=1}^{\infty} \langle C_{\alpha_1}, T_n C_{\alpha_2} \rangle q^n, \tag{2.5}$$

*and this is an elliptic modular form of level $\Gamma_0(ND_B)$.*

*Proof*: The first part is immediate from the definition of $\Theta_{\mathrm{KM}}$ and Lemma 2.8. To determine the level, it suffices to show that $\mathbb{1}_{\hat{\mathcal{O}}}$ is invariant under $\Gamma_0(ND_B)$ under the Weil representation. This computation proceeds exactly as in [Bra25, Proposition 3.1], since the dual lattice of $\mathcal{O}$ with respect to the bilinear form $\langle b_1, b_2 \rangle = \mathrm{Tr}\bigl(b_1 \overline{b_2}\bigr)$ is $\frac{1}{D_B N} \cdot \mathcal{O}$. □

This gives a new proof of the modularity of Richards' generating series.

## 3 A real quadratic theta lift

Having dealt with one half of the seesaw diagram in Figure 1, we now turn to the other.

### 3.1 An $F$-quadratic form on $B$

Recall that we fixed a pair of embeddings $\alpha_j : F_j \to B$, $j \in \{1, 2\}$. Consider the involution $\varepsilon$ on $L := F_1 \otimes F_2$ which restricts to the Galois involutions on the respective factors:

$$\varepsilon\bigl(\sqrt{D_1} \otimes 1\bigr) = -\sqrt{D_1} \otimes 1 \quad \text{and} \quad \varepsilon\bigl(1 \otimes \sqrt{D_2}\bigr) = -1 \otimes \sqrt{D_2}. \tag{3.1}$$

The fixed algebra $F = (F_1 \otimes F_2)^{\varepsilon=1}$ fits into the following diagram of étale algebras:

$$\begin{array}{c} L \\ {\diagup}\ |\ {\diagdown} \\ F_1 \quad F \quad F_2 \\ {\diagdown}\ |\ {\diagup} \\ \mathbb{Q} \end{array}$$

If $F_1$ and $F_2$ are linearly disjoint, i.e. if their discriminants are coprime, then $L = F_1 F_2$ and $F$ is the third real quadratic field in this biquadratic extension of $\mathbb{Q}$. Note that the action $(x \otimes y) \cdot b = \alpha_1(x) b \alpha_2(y)$ defines an $L$-module structure on $B$, and in this way we view $B$ as an $F$-module.

**Remark 3.1**: *In the case where $B = \mathrm{Mat}_2(\mathbb{Q})$, which is excluded by hypothesis, we may also take $F_1 = \mathbb{Q} \oplus \mathbb{Q}$, the "split" real quadratic field. If $F \cong F_2$ is a real quadratic field and $\mathcal{O}$ is the standard Eichler order of level $p$, the natural analogue of Theorem 1.2 is proved in [Bra23a]. This case is more complicated because the theta integral is no longer absolutely convergent, and requires a regularisation argument.*

**Proposition 3.2**: *The inner product $\langle b_1, b_2 \rangle = 2\,\mathrm{Tr}\!\left(b_1 \overline{b_2}\right)$ on $B$ may be lifted to a unique $F$-bilinear form $\langle b_1, b_2 \rangle_F$ satisfying*

$$\mathrm{Tr}_{F/\mathbb{Q}} \langle b_1, b_2 \rangle_F = \langle b_1, b_2 \rangle. \tag{3.2}$$

*Proof*: Existence and uniqueness can be proved abstractly using the bijectivity of the trace map $\mathrm{Tr}_{F/\mathbb{Q}} : \mathrm{Hom}_F(B, F) \to \mathrm{Hom}_{\mathbb{Q}}(B, \mathbb{Q})$ as in [Phi15, Lemma 3.1.5]. It is also straightforward to determine an explicit formula by writing $\langle b_1, b_2 \rangle_F = u + v\sqrt{D}$ and solving for $u$ and $v$; if $D = D_1 D_2$ this gives

$$\langle b_1, b_2 \rangle_F := \frac{1}{2} \langle b_1, b_2 \rangle + \frac{1}{2\sqrt{D}} \langle \sqrt{D} \cdot b_1, b_2 \rangle. \tag{3.3}$$

$\square$

Let $q_F(b) = \frac{1}{2} \langle b, b \rangle_F$ denote the associated $F$-quadratic form on $B$.

**Proposition 3.3**: *There exists an isometry $\iota_L : (B, q_F) \cong (L, q_\alpha)$, where $q_\alpha$ is the quadratic form $q_\alpha(x) = \alpha \, \mathrm{Nm}_{L/F}(x)$ and $\alpha = q_F(1)$.*

*Proof*: First note that

$$\langle x \cdot b_1, b_2 \rangle_F = \langle b_1, \varepsilon(x) \cdot b_2 \rangle_F, \tag{3.4}$$

for all $x \in L$ and $b_1, b_2 \in B$. The argument in the proof of [Mil69, Lemma 1.1] implies that $\langle -, - \rangle_F$ lifts to an $L$-hermitian form $h$ on $B$ such that $\mathrm{Tr}_{L/F} h = \langle -, - \rangle_F$, or equivalently, $h(b, b) = q_F(b)$ for all $b \in B$. By the classification of one-dimensional $L$-hermitian spaces, the vector space isomorphism $L \to B$ given by $x \mapsto x \cdot 1$ lifts to an isomorphism of hermitian spaces $(L, \alpha x \varepsilon(y)) \to (B, h)$ for some $\alpha \in F^\times$. By construction, $\alpha \cdot 1 = h(1, 1) = q_F(1)$.
$\square$

**Corollary 3.4**: *For any $\varepsilon_i \in \mathcal{O}_i^{\mathrm{Nm}=1}$ and $b \in B$, we have $q_F(\alpha_i(\varepsilon_i) \cdot b) = q_F(b)$.*

**Remark 3.5**:

(i) *Naively, one might hope that the form $q_F$ only depends on each embedding $\alpha_j$ up to equivalence of optimal embeddings, i.e. conjugation by elements of $\mathcal{O}^1$. However, it is not hard to construct examples numerically showing that this is false. On the other hand, $\alpha$ is invariant under simultaneous conjugation of $\alpha_1$ and $\alpha_2$.*

(ii) *Since $\iota_L(\mathcal{O})$ is a rank 4 lattice preserved by multiplication by $\mathcal{O}_1 \otimes \mathcal{O}_2$, it is a fractional ideal.*

(iii) *One can give an alternative construction of $q_F$ in terms of the modular cross ratio, see [Daa25, Proposition 4.4] for an argument in the imaginary quadratic setting.*

The isometry in Proposition 3.3 gives an isomorphism $\mathrm{SO}_{(B,q_F)} \xrightarrow{\sim} \mathrm{SO}_{(L,q_\alpha)}$ of algebraic groups over $F$. For convenience, we set $S := \mathrm{SO}_{(L,q_\alpha)}$.

**Proposition 3.6**: *The map $\eta : (t_1, t_2) \mapsto t_1 \otimes t_2^{-1}$ determines a surjection of groups $T(\mathbb{Q}) \cong S(F)$ with kernel $\mathbb{Q}^\times$.*

*Proof*: It is well-known that $\mathrm{SO}_{(L, \mathrm{Nm}_{L/F})} \cong L^1 := \{x \in L : x \cdot \varepsilon(x) = 1\}$, and rescaling $\mathrm{Nm}_{L/F}$ by $\alpha$ does not alter the situation. It therefore suffices to show that $T(\mathbb{Q})$ fits into the short exact sequence

$$1 \to \mathbb{Q}^\times \to T(\mathbb{Q}) \to L^1 \to 1. \tag{3.5}$$

If $t = (t_1, t_2) \in T(\mathbb{Q})$, then

$$\eta(t) \cdot \varepsilon(\eta(t)) = \frac{\mathrm{Nm}(t_1)}{\mathrm{Nm}(t_2)} = 1, \tag{3.6}$$

so the map $T(\mathbb{Q}) \to L^1$ is well-defined, and has kernel equal to $\mathbb{Q}^\times$. Surjectivity follows from the argument in [Bra25, Proposition 2.2]. □

In fact, one can show that the sequence

$$1 \to k^\times \to T(k) \to S(k \otimes F) \to 1 \tag{3.7}$$

is exact when $k$ is either a characteristic 0 field, $\mathbb{A}$ or $\mathbb{A}^\infty$. The proof is identical to that of [HY12, Proposition 2.13].

Since $\eta$ is the restriction of the natural map

$$\mathrm{GSpin}_B(\mathbb{Q}) \to \mathrm{SO}_B(\mathbb{Q}) \cong \mathrm{SO}_L(\mathbb{Q}), \tag{3.8}$$

we obtain the following diagram of groups with exact rows:

$$\begin{array}{ccccccccc}
1 & \longrightarrow & \mathbb{Q}^\times & \longrightarrow & \mathrm{GSpin}_B(\mathbb{Q}) & \longrightarrow & \mathrm{SO}_L(\mathbb{Q}) & \longrightarrow & 1 \\
& & \| & & \uparrow & & \uparrow & & \\
1 & \longrightarrow & \mathbb{Q}^\times & \longrightarrow & T(\mathbb{Q}) & \xrightarrow{\eta} & S(F) & \longrightarrow & 1
\end{array}$$

Let $X_S := S(F) \backslash (S(\mathbb{A}_F^\infty) \times \mathbb{D}_S^+/K_S^\infty)$, where $K_S^\infty \subset S(\mathbb{A}_F^\infty)$ is the image of $K_T^\infty$ under map induced by $\eta$.

**Lemma 3.7**: *We have an isomorphism $X_T \cong X_S$, and in particular*

$$X_S = \bigsqcup_{h \in \mathrm{Cl}^+ T} \Gamma_S \backslash \mathbb{D}_S^+, \tag{3.9}$$

*where $\Gamma_S := \eta(\Gamma_T)$.*

*Proof*: This is a straightforward consequence of Equation (3.7). □

## 3.2 The Kudla–Millson form restricted to $S$

The pair $\left(S, \mathrm{Res}_{F/\mathbb{Q}} \mathrm{SL}_2\right)$ forms a weakly dual reductive pair in the terminology of [Kud84]. The oscillator representation on the space $\mathcal{S}(L \otimes_F \mathbb{A}_F)$ extends the action of $\omega$ under the embedding $\mathrm{SL}_2(\mathbb{A}) \to \mathrm{SL}_2(\mathbb{A}_F)$, so we denote it by the same symbol. Similarly, the local representations on $\mathcal{S}(L \otimes_F F_v)$ for any place $v$ of $F$ will also be denoted by $\omega$, as before.

By functoriality proved in [Bra23a, §5], the Kudla–Millson form on $\mathbb{D}_S$ is the pullback of that on $\mathbb{D}$ by the inclusion $\mathbb{D}_S \to \mathbb{D}$. Writing $F_\infty := F \otimes_\mathbb{Q} \mathbb{R} = F_\sigma \times F_{\sigma'}$ where $\sigma, \sigma' : F \hookrightarrow \mathbb{R}$, we get a natural decomposition $L \otimes_\mathbb{Q} \mathbb{R} \cong L_\sigma \times L_{\sigma'}$ of quadratic spaces over $\mathbb{R}$. Here $L_\sigma$ is the quadratic space over $F_\sigma \cong \mathbb{R}$ with quadratic form $q_\sigma(x) = \sigma(\alpha) \mathrm{Nm}_{L_\sigma/F_\sigma}(x)$, of signature $(1,1)$.

Let $\mathbb{R}^{1,1}$ denote the 2-dimensional quadratic space over $\mathbb{R}$ with quadratic form $q(x,y) := x \cdot y$, and notice that $\mathbb{R}^\times \cong \mathrm{SO}_{\mathbb{R}^{1,1}}$ via

$$t \mapsto \begin{pmatrix} t & 0 \\ 0 & t^{-1} \end{pmatrix}. \tag{3.1}$$

**Lemma 3.8**: *For each embedding $\sigma : F \to \mathbb{R}$, the map*

$$\iota_\sigma : (x,y) \mapsto \left(\sqrt{|\sigma(\alpha)|} \cdot x, \mathrm{sgn}(\sigma(\alpha)) \cdot \sqrt{|\sigma(\alpha)|} \cdot y\right) \tag{3.2}$$

*defines an isometry $L_\sigma \to \mathbb{R}^{1,1}$.*

The Kudla–Millson Schwartz form for $\mathbb{R}^{1,1}$ was computed explicitly in [Bra23b]:

**Lemma 3.9** ([Bra23b], §5): *The Kudla–Millson Schwartz form $\varphi_{\mathrm{KM}}^{1,1}$ on the symmetric space of $\mathrm{SO}_{\mathbb{R}^{1,1}}$ is given by*[3]

$$\varphi_{\mathrm{KM}}^{1,1}(x,y) := \left(\frac{x}{t} + yt\right) \cdot e^{-\pi\left(\left(\frac{x}{t}\right)^2 + (yt)^2\right)} \frac{\mathrm{d}t}{t}. \tag{3.3}$$

The symmetric space $\mathbb{D}_\sigma^+$ is isomorphic to $\mathbb{R}_{>0}$, and so

$$\mathbb{D}_S^+ \cong \mathbb{D}_\sigma^+ \times \mathbb{D}_{\sigma'}^+ \cong \mathbb{R}_{>0} \times \mathbb{R}_{>0}, \tag{3.4}$$

where the last map is induced by $\iota_\sigma \times \iota_{\sigma'}$. Let $t_\sigma$ be the coordinate on $\mathbb{D}_\sigma^+$, and $\mathrm{d}t_\sigma/t_\sigma$ the invariant differential. Note that $\mathrm{SO}_{L_\sigma} \cong F_\sigma^\times \cong \mathbb{R}^\times$ via the isomorphism

$$\begin{pmatrix} t_\sigma & 0 \\ 0 & t_\sigma^{-1} \end{pmatrix} \mapsto t_\sigma. \tag{3.5}$$

If $\varphi_{\mathrm{KM}}^\sigma$ denotes the Kudla–Millson form on $\mathbb{D}_\sigma^+$, then it follows from functoriality of Branchereau's construction that $\varphi_{\mathrm{KM}}^\sigma = \iota_\sigma^* \varphi_{\mathrm{KM}}^{1,1}$, hence[4]

$$\varphi_{\mathrm{KM}}^\sigma(x,y) := \sqrt{|\sigma(\alpha)|}\left(\frac{x}{t_\sigma} + \mathrm{sgn}(\sigma(\alpha))yt_\sigma\right) e^{-\pi|\sigma(\alpha)|\left(\left(\frac{x}{t_\sigma}\right)^2 + (yt_\sigma)^2\right)} \frac{\mathrm{d}t_\sigma}{t_\sigma}. \tag{3.6}$$

It is convenient to rewrite this in terms of oscillator representation on $\mathcal{S}(L \otimes_F F_\sigma)$; to this end, define the Schwartz function

$$\varphi_\infty(x,y) := (x+y)e^{-\pi(x^2+y^2)}, \tag{3.7}$$

and set $\varphi_\sigma := \iota_\sigma^* \varphi_\infty$.

**Lemma 3.10**: *The the restriction of $\varphi_{\mathrm{KM}}$ to $\mathbb{D}_S^+$ is given by*

$$\varphi_{\mathrm{KM}}(x_\sigma, x_{\sigma'}) = \left(\omega(t_\sigma)\varphi_\sigma(x_\sigma)\frac{\mathrm{d}t_\sigma}{t_\sigma}\right) \wedge \left(\omega(t_{\sigma'})\varphi_{\sigma'}(x_{\sigma'})\frac{\mathrm{d}t_{\sigma'}}{t_{\sigma'}}\right). \tag{3.8}$$

*Proof*: By [Bra23b, Equation (5.10)], the restriction of $\varphi_{\mathrm{KM}}$ to $\mathbb{D}_S^+$ is given by the wedge product of the corresponding forms $\varphi_{\mathrm{KM}}^\sigma$ and $\varphi_{\mathrm{KM}}^{\sigma'}$. Now Equation (3.6) implies that $\omega(t_\sigma)\varphi_\sigma(x_\sigma) = \varphi_\sigma(t_\sigma^{-1}x, t_\sigma y)$ when $x_\sigma = (x,y)$. □

We now turn to the finite places. Recall that we defined orders $\mathcal{O}_j := \alpha_j^{-1}(\mathcal{O}) \subset F_j$. Via the isomorphism $L \to B$, the Eichler order $\mathcal{O}$ is naturally identified with $\mathcal{O}_1 \otimes \mathcal{O}_2$ as a lattice. Since this is stable under multiplication by $\mathcal{O}_{12} := \mathcal{O}_1 \otimes \mathcal{O}_2 \cap F$, the characteristic function $\varphi^\infty = \mathbb{1}_{\hat{\mathcal{O}}} \in \mathcal{S}(B_{\mathbb{A}^\infty})$ can be identified with a Schwartz–Bruhat function in $\mathcal{S}(L \otimes \mathbb{A}_F^\infty)$.

---

[3]There is a sign error in the calculation in [Bra23b], so his form should strictly speaking be $(-1) \cdot \varphi_{\mathrm{KM}}^{1,1}$. A short computation with Howe operators shows that our convention coincides with the Kudla–Millson form as constructed by Kudla and Millson. As we eventually take the wedge product of two such forms, the difference in signs is unimportant here.

[4]Our formula is slightly more complicated than that in [Bra25, Equation (3.78)] since $\alpha$ might not be totally positive — the results there should also be modified accordingly.

For $\tau_\sigma = x_\sigma + iy_\sigma \in \mathfrak{h}$ and $\underline{\tau} = (\tau_\sigma)_{\sigma:F\to\mathbb{R}}$, we may define a Kudla–Millson theta kernel over $\mathbb{A}_F$ by

$$\Theta_{\mathrm{KM}}^S(\underline{\tau}) := \sum_{v \in L} \omega(g_{\underline{\tau}})(\varphi^\infty \otimes \varphi_{\mathrm{KM}})(v), \tag{3.9}$$

where $g_{\underline{\tau}} \in \mathrm{SL}_2(\mathbb{A}_F)$ denotes the adelic matrix whose finite place components are 1, and infinite components are $(g_{\tau_\sigma})_\sigma$, with $g_{\tau_\sigma}$ the matrix in Equation (2.3). Let $i_\Delta : \mathrm{SL}_2(\mathbb{A}) \hookrightarrow \mathrm{SL}_2(\mathbb{A}_F)$ be natural embedding, and notice that $i_\Delta(1,...,g_\tau) = (1,...,g_\tau, g_\tau)$.

**Corollary 3.11**: *We have*

$$\int_{C_1 \times C_2} \Theta_{\mathrm{KM}}(g_\tau) = \int_{C_1 \times C_2} \Theta_{\mathrm{KM}}^S(\tau, \tau). \tag{3.10}$$

In the next section, we will show that the right-hand side of Equation (3.10) is the diagonal restriction of a Hilbert modular form.

### 3.3 Fourier coefficients

We next show that the Kudla–Millson theta kernel $\Theta_{\mathrm{KM}}^S$ gives rise to a Hilbert modular form of parallel weight 1, and compute it's Fourier coefficients explicitly. A good reference for Hilbert modular forms is [Gar90].

**Proposition 3.12**: *For any class $C \in H_2(X_T, \mathbb{Z})$, the function*

$$\vartheta_C(\underline{\tau}) := \int_C \Theta_{\mathrm{KM}}^S(\underline{\tau}) \tag{3.1}$$

*defines a Hilbert modular form of parallel weight* 1.

*Proof*: Since $\vartheta_C$ is a Kudla–Millson theta lift over $F$, this follows from [KM90, Theorem 2(bis.)]. □

Taking $C = C_{\alpha_1} \times C_{\alpha_2}$, Corollary 3.11 and Proposition 2.9 imply:

**Corollary 3.13** (Theorem 1.2): *For any $\tau \in \mathfrak{h}$, we have*

$$\vartheta := \vartheta_{C_{\alpha_1} \times C_{\alpha_2}}(\iota_\Delta(g_\tau)) = \sum_{n=1}^\infty \langle C_{\alpha_1}, T_n C_{\alpha_2} \rangle q^n. \tag{3.2}$$

The Fourier series of $\vartheta$ can be computed in a relatively straightforward manner:

**Proposition 3.14**: *We have*

$$\vartheta(\underline{\tau}) = \sum_{b \in \mathcal{O}_1^1 \backslash \mathcal{O}/\mathcal{O}_2^1} e^{2\pi i \operatorname{Tr}(q_F(b) \cdot \underline{\tau})} \cdot \mathbb{1}_{q_F(b) \gg 0} \cdot \varsigma(b), \tag{3.3}$$

*where* $\varsigma(b) \in \{\pm 1\}$ *is a sign depending on* $b$; *if* $q_F(1) \gg 0$, *then* $\varsigma(b) = \operatorname{Nm}_{F/\mathbb{Q}} \operatorname{Tr}_{L/F}(\iota(b))$.

**Remark 3.15**: *It is natural to expect* $\varsigma(b)$ *to be closely related to the sign of the intersection of the geodesics* $C_1$ *and* $b \cdot C_2$, *as defined in [Ric21, §1], but we were unable to find a satisfying proof.*

We need the following well-known lemma:

**Lemma 3.16**: *For any* $x, y \in \mathbb{R}$ *and* $\varphi_\infty(x, y) = (x + y)e^{-\pi(x^2 + y^2)}$, *we have*

$$\int_0^\infty \varphi_\infty(t^{-1}x, ty) \frac{\mathrm{d}t}{t} = \begin{cases} \operatorname{sgn}(x)e^{-2\pi xy} & \text{if } xy > 0, \\ 0 & \text{otherwise}. \end{cases} \tag{3.4}$$

For lack of a good reference, we give a quick proof:

*Proof*: Let $I(x, y)$ denote the integral in Equation (3.4). The change of variables $t = \sqrt{r \cdot |x/y|}$ gives

$$\begin{aligned} I(x, y) &= \int_0^\infty \varphi_\infty\left(\operatorname{sgn}(x)\sqrt{|xy|}r^{-\frac{1}{2}}, \operatorname{sgn}(y)\sqrt{|xy|}r^{\frac{1}{2}}\right) \frac{\mathrm{d}r}{2r} \\ &= \frac{\sqrt{|xy|}}{2} \int_0^\infty \left(\operatorname{sgn}(x)r^{-\frac{1}{2}} + \operatorname{sgn}(y)r^{\frac{1}{2}}\right) e^{-\pi|xy|(r+\frac{1}{r})} \frac{\mathrm{d}r}{r}. \end{aligned} \tag{3.5}$$

The modified Bessel function of the second kind, $B_\mu(z)$, has the integral representation

$$B_\mu(z) = \frac{1}{2} \int_0^\infty e^{-\frac{z}{2}(r+\frac{1}{r})} r^\mu \frac{\mathrm{d}r}{r}, \tag{3.6}$$

seen by setting $r = \exp(-t)$ in [Wat66, 6.22 (7)]. This gives

$$I(x, y) = \sqrt{|xy|} \cdot \left(\operatorname{sgn}(x) B_{-\frac{1}{2}}(2\pi|xy|) + \operatorname{sgn}(y) B_{\frac{1}{2}}(2\pi|xy|)\right). \tag{3.7}$$

Since $B_{-\mu}(z) = B_\mu(z)$, we see that $I(x, y) = 0$ if $\operatorname{sgn}(x) = -\operatorname{sgn}(y)$. Otherwise, applying the identity $B_{\frac{1}{2}}(z) = \sqrt{\frac{\pi}{2z}} e^{-z}$ from [Wat66, Equation 3.8.13], we find

$$I(x, y) = \operatorname{sgn}(x) e^{-2\pi|xy|} = \operatorname{sgn}(x) e^{-2\pi xy}, \tag{3.8}$$

proving the claim. □

*Proof of Proposition 3.14*: Unpacking definitions, we find that

$$\vartheta_C = \int_{\Gamma_S \backslash \mathbb{D}_S} \sum_{b \in \mathcal{O}} \omega(g_\tau) \varphi_{\mathrm{KM}}(b)$$
$$= \int_{\eta(\varepsilon_1^\mathbb{Z} \times \varepsilon_2^\mathbb{Z}) \backslash \mathbb{D}_S^+} \sum_{b \in \mathcal{O}} \omega(g_{\tau_\sigma}) \varphi_{\mathrm{KM}}(b) \qquad (3.9)$$

By splitting the sum into cosets for $\varepsilon_1^\mathbb{Z} \backslash \mathcal{O} / \varepsilon_2^\mathbb{Z}$, we can unfold the integral:

$$\vartheta_C = \int_{\mathbb{D}_S^+} \sum_{b \in \varepsilon_1^\mathbb{Z} \backslash \mathcal{O} / \varepsilon_2^\mathbb{Z}} \omega(g_\tau) \varphi_{\mathrm{KM}}(b)$$
$$= \sum_{b \in \varepsilon_1^\mathbb{Z} \backslash \mathcal{O} / \varepsilon_2^\mathbb{Z}} \prod_{\sigma: F \to \mathbb{R}} \int_{\mathbb{D}_\sigma^+} \omega(g_{\tau_\sigma}) \varphi_\sigma(t_\sigma^{-1} \cdot b) \frac{\mathrm{d}t_\sigma}{t_\sigma}, \qquad (3.10)$$

where we recall that $\iota_\sigma$ is the isometry $L_\sigma \xrightarrow{\sim} \mathbb{R}^{1,1}$ from Lemma 3.8.

Denote by $I_\sigma(b)$ the factor corresponding to $\sigma$. The formulas for the Weil representation give

$$\omega(g_{\tau_\sigma}) \varphi_\sigma(t_\sigma^{-1} \cdot b) = e^{2\pi i q_\sigma(b) x_\sigma} \varphi_\sigma(t_\sigma^{-1} \cdot \sqrt{y_\sigma} b), \qquad (3.11)$$

so

$$I_\sigma(b) = \int_{\mathbb{D}_\sigma^+} e^{2\pi i q_\sigma(b) x_\sigma} \cdot \varphi_\sigma(t_\sigma^{-1} \cdot \sqrt{y_\sigma} b) \frac{\mathrm{d}t_\sigma}{t_\sigma}. \qquad (3.12)$$

Now we use that $\varphi_\sigma = \iota_\sigma^* \varphi_\infty$, so

$$I_\sigma(b) = \int_{\mathbb{R}_{>0}} \varphi_\infty(t^{-1} \cdot \iota_\sigma(\sqrt{y_\sigma} b)) \frac{\mathrm{d}t}{t}. \qquad (3.13)$$

To apply Lemma 3.16, first let $(b_1, b_2)$ denote the image of $b$ in $L_\sigma \cong F_\sigma \times F_\sigma$ as If $(x, y) = \iota_\sigma(\sqrt{y_\sigma} b)$, then $\mathrm{sgn}(x) = \mathrm{sgn}(b_1)$, while $\mathrm{sgn}(y) = \mathrm{sgn}(\sigma(\alpha)) \mathrm{sgn}(b_2)$. Note that if $\mathrm{sgn}(\sigma(\alpha)) = 1$ and $q_\sigma(b) > 0$, then $\mathrm{sgn}(b_1)$ equals $\mathrm{sgn}(\mathrm{Tr}_{L_\sigma/F_\sigma} b)$. Let $\varsigma_\sigma(b) := \mathrm{sgn}(b_1)$, and $\varsigma(b) := \prod_{\sigma: F \to \mathbb{R}} \varsigma_\sigma(b)$.

Clearly $x \cdot y = q_\sigma(b)$, so Equation (3.4) gives

$$I_\sigma(b) = e^{2\pi i q_\sigma(b) x_\sigma} e^{-2\pi q_\sigma(b) y_\sigma} \cdot \mathbb{1}_{q_\sigma(b) > 0} \cdot \varsigma_\sigma(b)$$
$$= e^{2\pi i q_\sigma(b) \tau_\sigma} \cdot \mathbb{1}_{q_\sigma(b) > 0} \cdot \varsigma_\sigma(b). \qquad (3.14)$$

Thus

$$\vartheta(\underline{\tau}) = \sum_{b \in \varepsilon_1^\mathbb{Z} \backslash \mathcal{O} / \varepsilon_2^\mathbb{Z}} \prod_{\sigma: F \to \mathbb{R}} I_\sigma(b) = \sum_{b \in \varepsilon_1^\mathbb{Z} \backslash \mathcal{O} / \varepsilon_2^\mathbb{Z}} e^{2\pi i \, \mathrm{Tr}(q_F(b) \cdot \underline{\tau})} \cdot \mathbb{1}_{q_F(b) \gg 0} \cdot \varsigma(b), \qquad (3.15)$$

as claimed. □

Taking the diagonal restriction and comparing Fourier coefficients, we find:

**Corollary 3.17**: *We have*
$$\langle C_{\alpha_1}, T_n C_{\alpha_2}\rangle = \sum_{\substack{b\in \varepsilon_1^{\mathbb{Z}}\backslash \mathcal{O}/\varepsilon_2^{\mathbb{Z}} \\ q_F(b)\gg 0 \\ \operatorname{Tr}_{F/\mathbb{Q}} q_F(b)=n}} \varsigma(b). \tag{3.16}$$

When $\mathcal{O}_{12}$ is maximal, this reduces to Corollary 1.3.

# Bibliography

[Bra23a]   R. Branchereau, "Diagonal Restriction of Eisenstein Series and Kudla–Millson Theta Lift," *Forum Mathematicum*, vol. 35, no. 5, pp. 1373–1418, Sept. 2023, doi: 10.1515/forum-2022-0344.

[Bra23b]   R. Branchereau, "The Kudla–Millson Form via the Mathai–Quillen Formalism," *Canadian Journal of Mathematics*, pp. 1–26, Oct. 2023, doi: 10.4153/S0008414X23000573.

[Bra25]    R. Branchereau, "Kudla Millson Lift of Toric Cycles and Restriction of Hilbert Modular Forms," *Mathematische Zeitschrift*, vol. 309, no. 4, 2025, doi: 10.1007/s00209-025-03683-0.

[Daa25]    M. A. Daas, "CM-values of p-adic Θ-functions", *Research in the Mathematical Sciences*, vol. 12, no. 3, p. 53, Aug. 2025, doi: 10.1007/s40687-025-00521-x.

[Dam24]    H. Damm-Johnsen, "Diagonal Restrictions of Hilbert Modular Forms," DPhil Thesis, DPhil Thesis, University of Oxford, 2024.

[DPV21]    H. Darmon, A. Pozzi, and J. Vonk, "Diagonal restrictions of *p*-adic Eisenstein families", *Mathematische Annalen*, vol. 379, no. 1, pp. 503–548, Feb. 2021, doi: 10.1007/s00208-020-02086-2.

[DPV23]    H. Darmon, A. Pozzi, and J. Vonk, "The Values of the Dedekind–Rademacher Cocycle at Real Multiplication Points," *Journal of the European Mathematical Society*, May 2023, doi: 10.4171/jems/1344.

[DY13]     T. Du and T. Yang, "Quaternions and Kudla's Matching Principle," *Mathematical Research Letters*, vol. 20, no. 2, pp. 367–383, 2013, doi: 10.4310/MRL.2013.v20.n2.a12.

[Gar90]    P. B. Garrett, *Holomorphic Hilbert Modular Forms*. in The Wadsworth & Brooks/Cole Mathematics Series. Wadsworth & Brooks/Cole Advanced Books & Software, Pacific Grove, CA, 1990.

[GH81]     M. J. Greenberg and J. R. Harper, *Algebraic Topology - a First Course*. in Mathematical Lecture Note Series. Addison-Wesley, 1981.

[Hec17]    E. Hecke, "Über Die Kroneckersche Grenzformel Für Reelle Quadratische Körper Und Die Klassenzahl Relativ-Abelscher Körper," *Verhandl. d. Naturforschenden Gesell. i. Basel*, vol. 28, 1917.

[HY12]     B. Howard and T. Yang, "Singular Moduli Refined." Accessed: Aug. 29, 2023. [Online]. Available: http://arxiv.org/abs/1202.6410

[KM90]   S. S. Kudla and J. J. Millson, "Intersection Numbers of Cycles on Locally Symmetric Spaces and Fourier Coefficients of Holomorphic Modular Forms in Several Complex Variables," *Publications Mathématiques de l'IHÉS*, vol. 71, pp. 121–172, 1990.

[Kud03]  S. S. Kudla, "Integrals of Borcherds Forms," *Compositio Mathematica*, vol. 137, no. 3, pp. 293–349, July 2003, doi: 10.1023/A:1024127100993.

[Kud84]  S. Kudla, "Seesaw Dual Reductive Pairs, In: Automorphic Forms of Several Variables, Katata, 1983," *Progress in Mathematics*, vol. 46, pp. 244–268, 1984.

[Mil69]  J. Milnor, "On Isometries of Inner Product Spaces," *Inventiones mathematicae*, vol. 8, no. 2, pp. 83–97, June 1969, doi: 10.1007/BF01404612.

[Phi15]  A. Phillips, "Moduli of CM False Elliptic Curves," phdthesis, PhD Thesis, Boston College, 2015.

[Ric20]  J. Rickards, "Intersections of Closed Geodesics on Shimura Curves," phdthesis, PhD Thesis, McGill University, 2020.

[Ric21]  J. Rickards, "Computing Intersections of Closed Geodesics on the Modular Curve," *Journal of Number Theory*, vol. 225, pp. 374–408, Aug. 2021, doi: 10.1016/j.jnt.2020.11.024.

[Ric22]  J. Rickards, "Hecke Operators Acting on Optimal Embeddings in Indefinite Quaternion Algebras," *Acta Arithmetica*, vol. 204, no. 4, pp. 347–367, 2022.

[Shi75]  T. Shintani, "On Construction of Holomorphic Cusp Forms of Half Integral Weight," *Nagoya Mathematical Journal*, vol. 58, pp. 83–126, Sept. 1975, doi: 10.1017/S0027763000016706.

[Tu23]   L. W. Tu, "Lefschetz Fixed Point Theorems for Correspondences," *Mathematics Going Forward : Collected Mathematical Brushstrokes*. Springer International Publishing, Cham, pp. 47–57, 2023. doi: 10.1007/978-3-031-12244-6_5.

[Wat66]  G. N. Watson, *A Treatise on the Theory of Bessel Functions*, 2nd ed. in Cambridge Mathematical Library. Cambridge University Press, 1966, p. viii+804.

[Zha01]  S. Zhang, "Heights of Heegner Points on Shimura Curves," *Annals of Mathematics*, vol. 153, no. 1, pp. 27–147, 2001, doi: 10.2307/2661372.

*Email address:* havard-dj@proton.me